\documentclass[11pt,twoside]{article}
\usepackage{graphicx}
\usepackage{amssymb}
\usepackage{amsmath}
\usepackage{parskip}
\usepackage{url}
\usepackage{amscd}
\usepackage{graphicx}
\usepackage{float} 
\usepackage{cite}
\usepackage{amsmath}
\usepackage{enumitem}
\newtheorem{thm}{\noindent\bf Theorem}[section]
\newtheorem{pro}[thm]{\noindent\bf Proposition}
\newtheorem{qn}[thm]{\noindent\bf Question}
\newtheorem{lem}[thm]{\noindent\bf Lemma}
\newtheorem{clm}[thm]{\noindent\bf Claim}

\usepackage{graphicx}
\usepackage{rotating}
\usepackage{subcaption}
\usepackage{tikz}
\usepackage{pgfplots}
\pgfplotsset{compat=newest}
\usetikzlibrary{shapes.geometric}

\usetikzlibrary{calc}
\usetikzlibrary{decorations.pathreplacing}
\usetikzlibrary{positioning,patterns}
\usetikzlibrary{arrows,shapes,positioning}
\usetikzlibrary{decorations.markings}

\def \smvx {circle[radius = .07][fill = black]}

\tikzstyle{edge}=[very thick]
\definecolor{bostonuniversityred}{rgb}{0.8, 0.0, 0.0}
\definecolor{arsenic}{rgb}{0.23, 0.27, 0.29}
\tikzstyle{diredge}=[postaction={decorate,decoration={markings,
		mark=at position .65 with {\arrow[scale = 1]{stealth};}}}]

\tikzstyle{diredge2}=[postaction={decorate,decoration={markings,
		mark=at position .55 with {\arrow[scale = 1]{stealth};}}}]

\newcommand{\defPt}[3]{
	\def \pt {(#1, #2)}
	\coordinate [at = \pt, name = #3];
}

\tikzset{
   K5/.pic={
     \foreach \x in {1,...,5}{%
    \pgfmathparse{(\x-1)*360/5}
    \node[draw,circle,fill=black, inner sep=1 pt] (N-\x) at (\pgfmathresult:1 cm) [thick] {};
  }
  
  \foreach \x in {1,...,4}{%
    \foreach \y in {\x,...,5}{%
        \path (N-\x) edge[line width=1 pt] (N-\y);
  }
  }
  }
}

\tikzset{
   redK5/.pic={
     \foreach \x in {1,...,5}{%
    \pgfmathparse{(\x-1)*360/5}
    \node[draw,circle,fill=black, inner sep=1 pt] (N-\x) at (\pgfmathresult:1 cm) [thick] {};
  }
  
  \foreach \x in {1,...,4}{%
    \foreach \y in {\x,...,5}{%
        \path (N-\x) edge[line width=1 pt,red] (N-\y);
  }
  }
  }
}

\tikzset{
   blueK5/.pic={
     \foreach \x in {1,...,5}{%
    \pgfmathparse{(\x-1)*360/5}
    \node[draw,circle,fill=black, inner sep=1 pt] (N-\x) at (\pgfmathresult:1 cm) [thick] {};
  }
  
  \foreach \x in {1,...,4}{%
    \foreach \y in {\x,...,5}{%
        \path (N-\x) edge[line width=1 pt,blue] (N-\y);
  }
  }
  }
}

\tikzset{
   greenK5/.pic={
     \foreach \x in {1,...,5}{%
    \pgfmathparse{(\x-1)*360/5}
    \node[draw,circle,fill=black, inner sep=1 pt] (N-\x) at (\pgfmathresult:1 cm) [thick] {};
  }
  
  \foreach \x in {1,...,4}{%
    \foreach \y in {\x,...,5}{%
        \path (N-\x) edge[line width=1 pt,green] (N-\y);
  }
  }
  }
}

\newcommand{\fitellipsis}[3] 
{\draw [fill=white] let \p1=(#1), \p2=(#2), \n1={atan2(\y2-\y1,\x2-\x1)}, \n2={veclen(\y2-\y1,\x2-\x1)}
    in ($ (\p1)!0.5!(\p2) $) ellipse [ x radius=\n2/2+0.1cm, y radius=#3cm, rotate=\n1];
}

\newcommand{\fitellipsisnfill}[3] 
{\draw [] let \p1=(#1), \p2=(#2), \n1={atan2(\y2-\y1,\x2-\x1)}, \n2={veclen(\y2-\y1,\x2-\x1)}
    in ($ (\p1)!0.5!(\p2) $) ellipse [ x radius=\n2/2+0.1cm, y radius=#3cm, rotate=\n1];
}

\newcommand{\ceil}[1]{
    \left \lceil #1 \right \rceil
}
\newcommand{\floor}[1]{
    \left \lfloor #1 \right \rfloor
}

\begin{document}
\date{}
\title{{\Large\bf  The $(t-1)$-chromatic Ramsey number for paths }}

\author{{\normalsize Matija Buci\'{c} \quad Amir Khamseh}}

\maketitle \normalsize
\begin{abstract}
The following relaxation of the classical problem of determining Ramsey number of a fixed graph has first been proposed by Erd\H{o}s, Hajnal and Rado over 50 years ago. Given a graph $G$ and an integer $t \ge 2$ determine the minimum number $N$ such that in any $t$-coloured complete graph on $N$ vertices there is a copy of $G$ using only edges of some $t-1$ colours. We determine the answer precisely when $G$ is a path.
\end{abstract}
\maketitle

\section{Introduction}

Ramsey theory refers to a large body of results which all roughly speaking say that any sufficiently large structure must contain a well-organised substructure. For example, the classical theorem of Ramsey from 1929 \cite{ramsey} states that for any fixed graph $H$ in any $t$-edge colouring of a sufficiently large complete graph we can find a monochromatic copy of $H$. The $t$-colour Ramsey number of $H$, denoted $R^t(H)$, is defined as the smallest order of a complete graph satisfying this property.

Despite a large amount of attention over the years, there are still very few families of graphs for which the Ramsey numbers are known precisely. An early example of an exact Ramsey result, due to Gerencser and Gyarfas from 1967 \cite{gerencser}, determines $R^2(P_n)=\lfloor 3n/2-1 \rfloor,$ where $P_n$ denotes the path with $n$ vertices.   
Already the $3$-colour case was already substantially more difficult, \L{}uczak \cite{luczak99-con-match} determined $R^3(P_n)$ asymptotically, which was subsequently strengthened to an exact result for long enough paths by Gy\'arf\'as, Ruszink\'o, S\'ark\"{o}zy and Szemer\'edi \cite{GYRSS}. For four and more colours even the asymptotic result is still open, with best lower and upper bounds being due to Yongqi, Yuansheng, Feng and Bingxi \cite{Yongqi}, and Knierim and Su \cite{knierim}, respectively. 

Over the years, many generalisations of Ramsey numbers have been considered (the survey \cite{conlon2015recent} contains many examples); In particular, determining the generalised Ramsey number of a path in a variety of these more general settings has often turned out to be a very interesting problem \cite{alon,bucic1,bucic2,bucic3,bucic4,gyarfas4,hook2,hook,sudakov1}. A particularly natural generalisation first considered by Erd\H{o}s, Hajnal and Rado in 1965 \cite{EHR} asks the following. Given a graph $H$ and integers $ t\ge k \ge 1$ determine the minimum number $N$ such that in any $t$-coloured complete graph on $N$ vertices there is a copy of $H$ whose edges are coloured using at most $k$ different colours. Let us denote the answer to this question by $R^t_k(H)$. Observe that the case $k=1$ recovers the classical Ramsey problem and represents the most difficult instance of the problem. On the other hand, the case $k=t-1$ is the first non-trivial one and turns out to already be interesting.
In particular, $R^t_{t-1}$ is referred to as the $t-1$-chromatic Ramsey number and represents a natural first step towards understanding the classical Ramsey numbers.

Erd\H{o}s, Hajnal and Rado only considered the $(t-1)$-chromatic Ramsey numbers for complete graphs. In this direction Erd\H{o}s and Szemeredi \cite{ES} showed the answer for the complete graph $K_n$ is exponential in $n$ and the main open problem is to determine how the exponent behaves depending on the number of colours. The best known general bounds here are due to Alon, Erd\H{o}s, Gunderson and Molloy \cite{alon-erdos}.

First to explicitly ask the question of determining the $(t-1)$-chromatic Ramsey number for non-complete graphs were Chung and Liu in 1978. This question has been considered for a number of families of graphs over the years \cite{beam,budden,chli1,chli2,gyarfas5,gyarfas3,harborth,jacobson, kh2,kh,khomidiforests,khomidi,khomidi2}. In this paper we are interested in the $(t-1)$-chromatic Ramsey number of a path. Since the case of $t=2$ corresponds to the classical $2$-colour Ramsey problem, it is solved completely by the result of Gerencser and Gyarfas \cite{gerencser}. For $3 \le t\le 5$ the value of $R^t_{t-1}(P_\ell)$ was determined in \cite{meenakshi,khomidiforests,kh}. In this paper we determine the answer precisely for any number of colours, answering a question raised in \cite{khomidiforests}. 

\begin{thm}\label{main}
For $\ell,t\geq 2$ we have: $$R_{t-1}^t(P_\ell )=\ell+\left\lfloor\frac{\ell-2}{2^t-2}\right\rfloor.$$ 
\end{thm}

An interesting aspect of this result is that the $(t-1)$-chromatic Ramsey number of an even path matches that of the matching of the same order, determined by Gy\'arf\'as, S\'{a}rk\"{o}zy and Selkow \cite{gyarfas}, so in particular it also recovers this result. 

We also determine the correct value in the so called ``asymmetric'' variant of the problem where we seek paths of different lengths depending on which colour they omit. This more general result plays a key role in our proof of Theorem \ref{main}.

%
%
%

\section{Preliminaries}
Let us first agree on some convention and terminology. 
If $G$ is a graph, $V$ will denote its vertex set and $E$ its edge
set. The number of vertices of $G$ is denoted by $|G|$. 
$C_i$ denotes the cycle with $i$ vertices. A $t$-colouring will always refer to a colouring of the edges of $G$ using $t$ colours, or in other words a partition of $E$
into $t$ classes. Typically we use $[t]=\{1,2,\ldots,t\}$ as the set of colours. 

Let us now extend the definition of {{$(t-1)$-chromatic Ramsey numbers}}, to the asymmetric case. The $(t-1)$-chromatic Ramsey
number, denoted by  $R_{t-1}^t(G_1, G_2, \ldots, G_t)$, is defined to be the least number $n$ such that in any $t$-colouring of the complete graph $K_n$, for some $i$ we can find a copy of $G_i$ whose edges do not use colour $i$. Since we will focus on the case when each $G_i$ is a path let us also define $p(\ell_1\ldots, \ell_t)=R_{t-1}^t(P_{\ell_1}, \ldots, P_{\ell_t}).$

Since in $2$-colours this precisely recovers the classical Ramsey numbers the following is an immediate corollary of the result of Gerencser and Gy\'arf\'as \cite{gerencser}.
\begin{thm}\label{R_1,2}
If $2\leq \ell_1\leq \ell_2$, then $p({\ell_1},{\ell_2})=\ell_2+\lfloor \ell_1/2\rfloor-1$. 
\end{thm}

The following proposition gives us, mostly trivial bounds which relate the $t-1$-chromatic Ramsey numbers of graphs with respect to different number of colours. It is a generalization of results from Chung and Liu \cite{chli2}. 
\begin{lem}\label{gen}
Let $2 \le \ell_1, \ldots,
\ell_t$. Then we have $p(\ell_1,\ldots,\ell_t)\leq
p(\ell_1,\ldots,\ell_{t-1})$ and  equality holds if
$\ell_t \geq p(\ell_1,\ldots,\ell_{t-1})$.
\end{lem}
{\it Proof.} Let $n=p(\ell_1,\ldots,\ell_{t-1})$ and
$c:E(K_n)\rightarrow [t]$ be a $t$-colouring of the edges of $K_n$.
Define a new $(t-1)$-colouring $c'$ obtained by merging the colours $t-1$ and $t$. I.e.\
$c'(e)=i$ if $c(e)=i$, when $1\leq i\leq t-2$,  and
$c'(e)={t-1}$ if $c(e)=t-1$ or $c(e)=t$. 
By definition of $n$, in $K_n$ coloured according to $c'$ we can find a copy of $P_{\ell_i}$ avoiding colour $i$ for some $i \le t-1$. This $P_i$ also avoids colour $i$ in colouring $c$.  Hence  
$p(\ell_1,\ldots,\ell_t)\leq n$, 
as required.

Now suppose that $\ell_t\geq p(\ell_1,\ldots,\ell_{t-1})=n$. By definition, there exists a colouring of $K_{n-1}$ with $t-1$ colours $\{1,\dots,t-1\}$ 
such that $K_{n-1}$ does not contain $P_{\ell_i}$ avoiding colour $i$, for all $i: 1\leq i\leq t-1$.
One can view this as a $t$-colouring (not using colour $t$) of $K_{n-1}$ and since $n-1 < \ell_t$ we are guaranteed there is no path of order $\ell_t$ so this colouring also shows
$$n \le p(\ell_1,\ldots,\ell_t)\leq p(\ell_1,\ldots,\ell_{t-1})=n,$$
completing the proof.$\hfill\square$ 

We note that the above lemma does not require the sequence $\ell_i$ to be sorted.

\section{{{Proof of the main result}}}
In this section we will prove the following result.
\begin{thm}\label{AAAA}
Let $t\ge 3$ and $2\leq \ell_1\leq \ell_2\leq \ldots\leq \ell_{t-1}\leq \ell_t$,  then  
$$ 
p({\ell_1}, {\ell_2},\ldots,{\ell_{t}})= \begin{cases}
        \left\lfloor\frac{\ell_1+2\ell_2+\ldots+2^{t-1} \ell_{t}-2}{2^{t}-2}\right\rfloor & \text{if }\quad \ell_t<p({\ell_1}, {\ell_2},\ldots,{\ell_{t-1}})\\
        p({\ell_1}, {\ell_2},\ldots,{\ell_{t-1}}) & \text{else.}
\end{cases}
$$
\end{thm}
Given this result Theorem \ref{main} follows immediately by choosing each $\ell_j=\ell$.

Before proceeding to the proof, we introduce some convention and notation. Until further notice we suppose that $2\leq \ell_1\leq \ell_2\leq \ldots\leq \ell_{t-1}\leq \ell_t$. For $2\leq k\leq t$, we define 
$$s(\ell_1,\ldots, \ell_k):=\left\lfloor\frac{\ell_1+2\ell_2+\ldots+2^{k-1}\ell_{k}-2}{2^{k}-2}\right\rfloor.$$

Observe first that if $\ell_t \ge p({\ell_1}, {\ell_2},\ldots,{\ell_{t-1}})$ Lemma \ref{gen} implies $p({\ell_1}, {\ell_2},\ldots,{\ell_{t}})=p({\ell_1}, {\ell_2},\ldots,{\ell_{t-1}})$ so we only need to prove the first case of Theorem \ref{AAAA}. Namely provided that $\ell_t < p({\ell_1}, {\ell_2},\ldots,{\ell_{t-1}})$ we need to show  $p({\ell_1}, {\ell_2},\ldots,{\ell_{t}})=s(\ell_1,\ldots, \ell_t)$.

We begin with the lower bound in order to give an illustration on how the optimal examples behave, this will also serve to provide some motivation for the upper bound argument.
\subsection{The lower bound}
The following proposition contains the remaining case for the lower bound.

\begin{pro}\label{lower}
If $\ell_t<p({\ell_1}, {\ell_2},\ldots,{\ell_{t-1}})$,  then 
 $p({\ell_1}, {\ell_2},\ldots,{\ell_{t}})\geq  s(\ell_1,\ldots, \ell_t)$. 
\end{pro}
\noindent{{\bf Proof.}} Let $s=s(\ell_1,\ldots, \ell_t)$.
We shall show  $p({\ell_1}, {\ell_2},\ldots,{\ell_{t}})> s-1$ by exhibiting a $t$-colouring of $K_{s-1}$ which does not contain a colour $j$-avoiding path of order $\ell_j$ for all $1\le j \le t$. 

Let us first observe that if $s<\ell_t$ then $K_{s-1}$ can not contain $P_{\ell_t}$ regardless of the colouring so we can choose to only use the first $t-1$ colours. In particular, since $s<\ell_t<p({\ell_1}, {\ell_2},\ldots,{\ell_{t-1}})$ we can, by definition of $p({\ell_1}, {\ell_2},\ldots,{\ell_{t-1}})$ find a desired colouring of $K_{s-1}.$ From now on we assume $s \ge \ell_t$.

As we will see the structure of our colouring is very natural, although choosing the parameters correctly turns out to be a bit more technical.

We split the vertices of our $K_{s-1}$ into $t$ parts $A_1,\ldots,A_t$ with sizes $|A_j|=a_j.$ We colour in colour $j$ all edges contained within a part $A_j$. For edges between parts $A_i$ and $A_j$ with $i<j$ we colour the edges in colour $j$. See Figure \ref{fig:1} for an illustration.

The partitioning condition requires us to choose the sizes $a_j$ so as to satisfy $a_1+a_2+\ldots +a_t=s-1$. Let us now establish requirements on the sizes $a_j$ which guarantee there is no path avoiding colour $j$ of order $\ell_j$. Consider first the subgraph of $K_{s-1}$ consisting of edges coloured with colours in $[t]\setminus \{j\}$. We can split it into three parts $A=A_1 \cup \ldots \cup A_{j-1}; B=A_j; C= A_{j+1} \cup \ldots \cup A_t$ and we know $A\cup C$ spans a clique, all edges between $B$ and $C$ exist, $B$ spans an independent set and there are no edges between $A$ and $B$ (see Figure \ref{fig:2} for an illustration). Let us consider the longest path in this graph. If the path uses no vertex from $A$ its order is at most $2|C|+1.$ Otherwise by removing any vertex belonging to $A$ we split the path into pieces all of which have at least one endpoint in $C$, so in total have order at most $2|C|$. In total in this case the order is at most $|A|+2|C|.$ Putting the cases together the longest path we can find has order $\max(1,|A|)+2|C|$.

\begin{figure}
\begin{minipage}[t]{0.47\textwidth}
\centering
\captionsetup{width=\textwidth}
\begin{tikzpicture}[scale=0.75]
\defPt{0}{0}{l}
\defPt{3}{0}{r}
\defPt{6}{0}{t}

\foreach \i in {-0.75,-0.25,0.25,0.75}
{
{
    \foreach \j in {-0.75,-0.25,0.25,0.75}
        \draw[blue, line width =1pt ] ($(l)+(0,\i)$) -- ($(r)+(0,\j)$);
}
}

\foreach \i in {-0.75,-0.25,0.25,0.75}
{
    \foreach \j in {-0.75,-0.25,0.25,0.75}
        \draw[green, line width =1pt ] ($(r)+(0,\i)$) -- ($(t)+(0,\j)$);
}

\foreach \i in {-0.5, 0,0.5}
{
    \foreach \j in {-0.5, 0,0.5}
        \draw[green, line width =1pt ] ($(l)+(\i,0)$) to[out=110,in=70] ($(t)+(\j,0)$);
}

\fitellipsis{$(l)+(0,1)$}{$(l)-(0,1)$}{0.8};
\fitellipsis{$(r)+(0,1)$}{$(r)-(0,1)$}{0.8};
\fitellipsis{$(t)+(0,1)$}{$(t)-(0,1)$}{0.8};

\pic[rotate=18,scale=0.5] at ($(l)$) {redK5};
\pic[rotate=18, scale=0.5] at ($(r)$) {blueK5};
\pic[rotate=18, scale=0.5] at ($(t)$) {greenK5};

\node[] at ($(l)+(-0.9,-1.1)$) {$A_1$};
\node[] at ($(r)+(-0.9,-1.1)$) {$A_2$};
\node[] at ($(t)+(-0.9,-1.1)$) {$A_3$};

\end{tikzpicture}
\caption[First figure]{The $3$ colour example.}
\label{fig:1}
\end{minipage}\hfill
\begin{minipage}[t]{0.52\textwidth}
\centering
\begin{tikzpicture}[scale=0.75]
\defPt{0}{0}{l}
\defPt{3}{0}{r}
\defPt{6}{0}{t}

\foreach \i in {-0.75,-0.25,0.25,0.75}
{
    \foreach \j in {-0.75,-0.25,0.25,0.75}
        \draw[line width =1pt ] ($(r)+(0,\i)$) -- ($(t)+(0,\j)$);
}

\foreach \i in {-0.5, 0,0.5}
{
    \foreach \j in {-0.5, 0,0.5}
        \draw[line width =1pt ] ($(l)+(\i,0)$) to[out=110,in=70] ($(t)+(\j,0)$);
}

\fitellipsis{$(l)+(0,1)$}{$(l)-(0,1)$}{0.8};
\fitellipsis{$(r)+(0,1)$}{$(r)-(0,1)$}{0.8};
\fitellipsis{$(t)+(0,1)$}{$(t)-(0,1)$}{0.8};

\pic[rotate=18,scale=0.5] at ($(l)$) {K5};
\pic[rotate=18, scale=0.5] at ($(t)$) {K5};

\node[] at ($(l)+(-0.9,-1.1)$) {$A$};
\node[] at ($(r)+(-0.9,-1.1)$) {$B$};
\node[] at ($(t)+(-0.9,-1.1)$) {$C$};

\end{tikzpicture}
\captionsetup{width=\textwidth}
\caption{Edges allowed in a colour $2$-avoiding path.
}
\label{fig:2}
\end{minipage}
\end{figure}

Using this observation the claim will follow provided we can find $a_j$'s satisfying the following collection of equalities and inequalities, in addition to $a_1 > 0$ and $a_j\geq 0$ for all $j$.
\begin{align*} 
 2a_2+2a_3+2a_4+\ldots +2a_{t-1}+2a_{t}+1&\leq \ell_1-1, \\
 a_1+2a_3+2a_4+\ldots +2a_{t-1}+2a_{t}&\leq \ell_{2}-1, \\
 a_1+a_2+2a_4+ \ldots +2a_{t-1}+2a_{t}&= \ell_{3}-1, \\
  \vdots & \\
 a_1+a_2+\ldots +a_{t-3}+2a_{t-1}+2a_{t}&= \ell_{t-2}-1, \\
 a_1+a_2+a_3+\ldots +a_{t-2}+2a_{t}&= \ell_{t-1}-1,   \\
 a_1+a_2+a_3+\ldots +a_{t-2}+a_{t-1}&= \ell_t-1.   \\
a_1+a_2+a_3+\ldots +a_{t-2}+a_{t-1}+a_t &=s-1, \\
 \end{align*}
We will choose $a_t=s-\ell_t\ge 0$ and define recursively $a_i=2a_{i+1}+\ell_{i+1}-\ell_i$ for $3 \le i \le t-1$ the remaining two sizes $a_1,a_2$ we will choose later to satisfy the first two and the last of the above expressions. Note that $a_i \ge 0$ for these $i$ also follows since $\ell_{i+1}\ge \ell_i$. 

Note now that by repeatedly substituting $a_i=2a_{i+1}+\ell_{i+1}-\ell_i$ and using $a_t=s-\ell_t$ in the last equality we get:
\begin{align*}
    S:=a_3+a_4+\ldots+a_t& =-\ell_3+\ell_4+3a_4+a_5+\ldots+a_t\\
                         & =-\ell_3-2\ell_4+3\ell_5+7a_5+a_6+\ldots+a_t\\
                         & \:\: \vdots \\
                         & = -\ell_3-2\ell_4-2^2\ell_5-\ldots-2^{t-4}\ell_{t-1} +(2^{t-3}-1)\ell_{t}+(2^{t-2}-1)a_t\\
                         &= -\ell_3-2\ell_4-2^2\ell_5-\ldots-2^{t-3}\ell_{t} +(2^{t-2}-1)s
\end{align*}
So we now have to choose $a_1, a_2$ so as to satisfy $2a_2+2S+1\le \ell_1-1$, $a_1+2S\le \ell_2-1 $ and $a_1+a_2+S=s-1$. We choose:
\begin{equation*}
a_1= \ceil{\frac{2s-2S+\ell_2-\ell_1-1}{3}} \quad\quad\quad\quad\quad\quad\quad a_2=\floor{\frac{s-S+\ell_1-\ell_2-2}{3}}
\end{equation*}

Observe that removing the ceil sign in $a_1$ and floor in $a_2$ their sum is precisely $s-S-1$ which is an integer. This in turn implies the same holds for $a_1$ and $a_2$ since one comes with floor and one with ceil signs, namely $a_1+a_2+S=s-1$.

Observe now that by definition of $s$ we know $$(2^t-2)s\leq -2+\ell_1+2\ell_2+\ldots+2^{t-1}\ell_t=-2+\ell_1+2\ell_2+(2^{t}-4)s-4S,$$
which after rearranging implies $4S+2s+2 \le \ell_1+2\ell_2.$

We now use this inequality to verify the remaining required inequalities. First let us verify that $a_2 \ge 0$ (which then immediately implies also $a_1 \ge 1$ since $\ell_2 \ge \ell_1$). This follows from the above inequality since it gives
\begin{align*}
4(s-S)& \ge 6s-\ell_1-2\ell_2+2 \ge 4(\ell_2-\ell_1+2), 
\end{align*}
where the final inequality is equivalent to $6s+3\ell_1 \ge 6\ell_2+6,$ which holds since $s \ge \ell_t \ge \ell_2$ and $\ell_1 \ge 2$.

Once again rearranging our inequality gives us 
$$\ell_2-1 \ge \frac{2s-2S-1-\ell_1+\ell_2}{3}+2S.$$
Since the LHS is an integer we may take the ceil of the RHS which implies $\ell_2-1 \ge a_1+2S,$ as desired. We obtain the final inequality after yet another rearrangement gives $$\ell_1-1 \ge 2\cdot \frac{s-S+\ell_1-\ell_2-2}{3}+1+2S\ge 2a_2+2S+1.$$

\vspace{-1.1cm}
\hfill $\Box$
%
\subsection{The upper bound}
Our proof of Theorem \ref{AAAA} is going to proceed by induction on $\ell_1+\ldots+\ell_t$. The following lemma will serve as the basis. Throughout the section we write $s=s({\ell_1}, {\ell_2},\ldots,{\ell_{t}}).$
\begin{lem}\label{lemi<3}
If  $\ell_1\leq 3$, then $p({\ell_1}, {\ell_2},\ldots,{\ell_{t}})\leq s$. 
\end{lem}
\noindent{{\bf Proof.}} 
By repeated use of Lemma \ref{gen} and an application of Theorem \ref{R_1,2} we get
$$p({\ell_1}, {\ell_2},\ldots,{\ell_{t}}) \leq p({\ell_1}, {\ell_2})=\ell_2+\lfloor \ell_1/2\rfloor-1=\ell_2 \leq s,$$
where the last equality follows since $2\leq \ell_1\leq 3$ and the final inequality follows from definition of $s$ and $2^{t-1}\ell_2+\ldots+2^2\ell_2+2\ell_2 \leq 2^{t-1}\ell_t+\ldots+2^2\ell_3+2\ell_2+\ell_1-2$.$\hfill \Box$ 



%
We are now ready to prove the remaining case of the upper bound of Theorem \ref{AAAA} and hence complete its proof.
\begin{thm}\label{upper}
If  $t\geq 3$ and $\ell_t<p({\ell_1}, {\ell_2},\ldots,{\ell_{t-1}})$, then 
 $p({\ell_1}, {\ell_2},\ldots,{\ell_{t}})\leq s$. 
\end{thm}
\noindent{{\bf Proof.}} The proof is by induction on both $t$ and $\ell_1+\ldots+\ell_t$ with base cases being taken care of by Theorem \ref{R_1,2} and Lemma \ref{lemi<3}. Let us suppose now that $t \ge 3$ and $\ell_i \ge 4$. We assume that Theorem \ref{AAAA} holds for $t-1$ paths of any lengths and for $t$ paths of lengths $2\le \ell_1'\le\ldots,\ell_{t}'$ provided  $\ell_1'+\ldots+\ell_t'< \ell_1+\ldots+\ell_t$. In particular, we are going to use the following immediate corollary of the inductive assumption.
\begin{equation}\label{eq:induction}
    p(\ell_1',\ldots, \ell_t') \le \max\{s(\ell_1',\ldots, \ell_t'), \ell_t'\}.
\end{equation}

Suppose that the edges of $G=K_{s}$  are coloured by
$1, \ldots, t$. We split the proof in a sequence of claims. 
\begin{clm}\label{A4}
 $\ell_t\leq s$. 
\end{clm}
\noindent{\it{Proof of Claim \ref{A4}.}} 
Let us first show that $p({\ell_1}, {\ell_2},\ldots,{\ell_{t-1}})\leq s({\ell_1}, {\ell_2},\ldots,{\ell_{t-1}}).$

If $t \ge 4$ then using Lemma \ref{gen}, 
$\ell_{t-1}\le \ell_t<p({\ell_1}, {\ell_2},\ldots,{\ell_{t-1}})\leq p({\ell_1}, {\ell_2},\ldots,{\ell_{t-2}})$ so the inductive assumption implies $p({\ell_1}, {\ell_2},\ldots,{\ell_{t-1}})\leq s({\ell_1}, {\ell_2},\ldots,{\ell_{t-1}}).$ Note that the same holds if $t=3$ as well, by Theorem \ref{R_1,2}.

Since
$\ell_t< p({\ell_1}, {\ell_2},\ldots,{\ell_{t-1}}) \leq s({\ell_1}, {\ell_2},\ldots,{\ell_{t-1}})$,  we have
$\ell_t(2^{t-1}-2)\leq \ell_1+2\ell_2+\ldots+2^{t-2} \ell_{t-1}-2$.
Adding $2^{t-1}\ell_t$ to both sides and dividing by $2^t-2$ shows 
$$\ell_t \le \frac{\ell_1+2\ell_2+\ldots+2^{t-2} \ell_{t-1}-2}{2^t-2}.$$
Since $\ell_t$ is an integer we may take the floor of the RHS, which gives our desired conclusion.
$\hfill \dashv$ 

The following claim will be a useful estimate when applying induction.
\begin{clm}\label{removal}
Let $k \ge 1$ be an integer and let $\ell_i'$ be such that for some $j$ we have $\ell_j'=\ell_j$ and for all $i \neq j$ we have $\ell_i'=\ell_i-2k \ge 2$. Then $p(\ell_1',\ldots,\ell_t') \le \max\{\ell_1',\ldots, \ell_t', s-k\}.$
\end{clm}
\noindent{\it{Proof of Claim \ref{removal}.}} We note that we are not assuming that the sequence $\ell_i'$ is sorted. Indeed let $2\le \ell_1''\le \ldots \le \ell_t''$ be the sorted permutation of $\ell_1',\ldots, \ell_t'$.  The induction hypothesis, in the form of \eqref{eq:induction}, implies $p(\ell_1',\ldots,\ell_t') \le \max\{\ell_t'',s(\ell_1'',\ldots, \ell_t'')\}.$ To complete the claim observe that  
\begin{align*}
    s(\ell_1'',\ldots,\ell_t'')& = \floor{\frac{\ell_1''+2\ell_2''+\ldots+2^{t-1}\ell_t''-2}{2^t-2}}\\
    &\le \floor{\frac{\ell_1-2k+2(\ell_2-2k)+\ldots+2^{t-2}(\ell_{t-1}-2k)+2^{t-1}\ell_t-2}{2^t-2}}\\
    & = \floor{\frac{\ell_1+2\ell_2+\ldots+2^{t-1}\ell_t-2}{2^t-2}}-k=s-k.
\end{align*}
In the inequality above we first took out the $t-1$ terms $-2k$, note that their contribution is maximal if they are matched with the smallest possible powers of $2$. The remaining terms are products of $1,2,\ldots, 2^{t-1}$ with respectively with corresponding term of a certain permutation of $\ell_1,\ldots, \ell_t$. Since by the rearrangement inequality such a product is maximised when both sequences are ordered in the same way we know the contribution is indeed $\ell_1+2\ell_2+\ldots+2^{t-1}\ell_t.$ $\hfill \dashv$

The following claim provides us with a slightly stronger estimate which will be useful in certain regimes. 
\begin{clm} \label{induction}
Let $2q+2 \le \ell_1\le \ldots \le \ell_t$ and assume $s=q+\ell_j-1$ then
$$p(\ell_1-2q,\ldots, \ell_{j-1}-2q, \ell_{j+1}-2q, \ldots, \ell_t-2q) \le \max(\ell_t-2q,\ell_j-1).$$
\end{clm}
\noindent{\it{Proof of Claim \ref{induction}.}} As in the previous claim we know that the desired expression is either at most $\ell_t-2q$ or at most the floor of
\begin{align*}
&\frac{\ell_1-2q+\ldots+2^{j-2}(\ell_{j-1}-2q)+2^{j-1}(\ell_{j+1}-2q)+\ldots+2^{t-2}(\ell_{t}-2q)-2}{2^{t-1}-2}=\\
&\frac{\ell_1+\ldots+2^{j-2}\ell_{j-1}+2^{j-1}\ell_{j+1}+\ldots+2^{t-2}\ell_{t}-(2^t-2)q-2}{2^{t-1}-2} \le \\
& \frac{\ell_1+\ldots+2^{j-2}\ell_{j-1}+2^{j-1}\ell_{j}+2^{j}\ell_{j+1}+\ldots+2^{t-1}\ell_{t}-2^{t-1}\ell_j-(2^t-2)q-2}{2^{t-1}-2} < \\
& \frac{(s+1)(2^t-2)-2^{t-1}\ell_j-(2^t-2)q}{2^{t-1}-2}=\frac{(2^t-2)\ell_j-2^{t-1}\ell_j}{2^{t-1}-2}=\ell_j \\
\end{align*}
where in the first inequality we added $2^{i-2}(\ell_{i}-\ell_j)\ge 0$ for all $i>j$ and added and subtracted $2^{j-1}\ell_j$ this contributed the term of $(2^{j-1}+2^{j-1}+2^j+\ldots+2^{t-2})\ell_j=2^{t-1}\ell_j$. In the second inequality we used the definition of $s$ and in the penultimate equality we used $s=q+\ell_j-1$.
Since for our bound we take the floor of the first expression and the inequality is strict the expression is at most $\ell_j-1,$ as claimed. $\hfill \dashv$

The following two claims contain the main part of our argument. 
\begin{clm}\label{lemcycle}
If $G$ contains a cycle of length $\ell_j-1$ avoiding colour $j$ then $G$ contains $P_{\ell_{j'}}$ avoiding colour $j'$ for some $j'\in[t]$. 
\end{clm}
\noindent{\it{Proof of Claim \ref{lemcycle}.}}
Let $C$ be the cycle of length $\ell_j-1$ given by the lemma and
let $Q$ be the graph induced by the $q=s-(\ell_j-1)$ vertices of $G$ not in $C$. If there is any edge not of colour $j$ between $C$ and $Q$ then we can use it to extend $C$ into a $P_{\ell_j}$ avoiding colour $j$ and are done. Otherwise, all of the edges joining the $\ell_j-1$ vertices of $C$ and
the $q$ vertices of $Q$ have colour $j$. We can always use these edges to build a path in colour $j$ of order $2\min(\ell_j-1,q)$ and actually by one longer if $\ell_j-1 \neq q.$ Observe first that if $\ell_j-1=q$ then we can find a Hamilton path of colour $j$, so of order at least $s\ge \ell_t$ which provides the desired $j'$ avoiding path for any $j' \neq j$. So we may assume $\ell_j-1 \neq q$ and that we have a colour $j$ path of order $2\min(\ell_j-1,q)+1$. 

\vspace{0.5cm}
\textbf{\noindent Case 1.} $\ell_j-1 < q.$ 

The colour $j$ path has order $2\ell_j-1 \ge \ell_1$ so it serves as the desired path, avoiding colour $1$, unless $j=1$. In this case it would also serve as the desired path avoiding colour $2$ unless $\ell_2 > 2(\ell_1-1)+1$. This implies $\ell_2 \ge 2(\ell_1-1)+2$ and allows us to use Claim \ref{removal} to show \begin{equation*}
p({\ell_1},{\ell_2-2(\ell_1-1)},{\ell_3-2(\ell_1-1)}, \ldots, {\ell_t-2(\ell_1-1)})\leq \max(\ell_1,\ell_t-2(\ell_1-1),s-(\ell_1-1))\le q,  
\end{equation*}
where the final inequality holds since by the case assumption we have $\ell_1 \le q$ and since $\ell_t-2(\ell_1-1)\le s-2(\ell_1-1)=q-(\ell_1-1)$.

\begin{figure}[t] 
\begin{center}
\begin{tikzpicture}[scale=0.75]
\defPt{0}{0.95}{l}
\defPt{3}{0}{r}

\defPt{2.65}{-1}{p1}
\defPt{3.35}{-1}{p2}
\defPt{2.65}{-0.5}{p3}
\defPt{3.35}{-0.5}{p4}
\defPt{2.65}{0}{p5}
\defPt{3.35}{0}{p6}
\defPt{0}{0.5}{p7}
\defPt{3}{0.5}{p8}
\defPt{0}{1}{p9}
\defPt{3}{1}{p10}
\defPt{0}{1.5}{p11}
\defPt{3}{1.5}{p12}

\fitellipsis{$(l)+(0,0.8)$}{$(l)-(0,0.8)$}{0.7};
\fitellipsis{$(r)+(0,1.7)$}{$(r)-(0,1.7)$}{0.9};


\node[] at ($(l)+(-0.5,-1.2)$) {$C$};
\node[] at ($(r)+(1,-1.6)$) {$Q$};

\draw [thick,decorate,decoration={brace,mirror,amplitude=5pt,raise=2pt},yshift=0pt]($(l)+(-0.7,0.6)$) -- ($(l)+(-0.7,-0.6)$) node [black,midway,xshift=-0.85cm] {\small $\ell_1-1$};

\draw [thick,decorate,decoration={brace,amplitude=5pt,raise=2pt},yshift=0pt]($(l)+(4,0.6)$) -- ($(l)+(4,-0.6)$) node [black,midway,xshift=0.85cm] {\small $\ell_1-1$};

\draw [thick,decorate,decoration={brace,amplitude=5pt,raise=2pt},yshift=0pt]($(l)+(4,-0.9)$) -- ($(l)+(4,-2.1)$) node [black,midway,xshift=1.4cm] {\small $\ell_i-2(\ell_1-1)$};

\foreach \i in {1,...,12}
{
\draw[] (p\i) \smvx;
}

\draw[] (p2) -- (p1) -- (p4) -- (p3) -- (p6) -- (p5) -- (p7) -- (p8) -- (p9) -- (p10) -- (p11) -- (p12);
\end{tikzpicture}
\caption{Illustration on how we construct the desired path. The path of length $\ell_i-2(\ell_1-1)$ avoiding some colour $i\neq 1$ is found by induction inside $Q$. We extend it using as depicted using the fact that all edges between $C$ and $Q$ have colour $1$.}\label{figure1}
\end{center}
\end{figure}
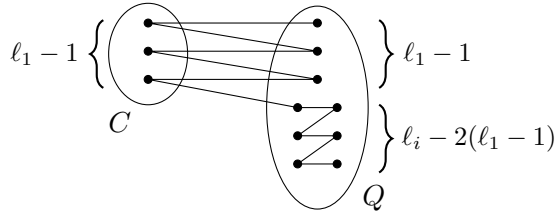

This means that we can either find the desired colour $1$ avoiding path of order $\ell_1$ inside $Q$ or we can find a path of order $\ell_i-2(\ell_1-1)$ avoiding some colour $i>1$ inside $Q$. Since  $\ell_i-2(\ell_1-1)\le s-2(\ell_1-1)=q-(\ell_1-1)$ there are $\ell_1-1$ vertices in $Q$ not on this path. This means that we can extend the path by going to $C$ and then alternating between these $\ell_1-1$ vertices of $Q$ not on the path and $C$ (see Figure \ref{figure1} for an illustration). This produces a path of order $\ell_i-2(\ell_1-1)+2(\ell_1-1)=\ell_i$ which avoids colour $i$ (recall that all edges between $C$ and $Q$ are in colour $j=1$) as desired. 

\vspace{.5cm} \noindent{\bf  Case $2$.} $\ell_j-1 > q.$

This means that in between $Q$ and $C$ we can find a colour $j$ path of order $2q+1$. This means that we are done unless $\ell_i \ge 2q+2$ for all $i \neq j$. This allows us to use Claim \ref{induction} which gives  
$$p(\ell_1-2q,\ldots, \ell_{j-1}-2q, \ell_{j+1}-2q, \ldots, \ell_t-2q) \le \max(\ell_t-2q,\ell_j-1)=\ell_j-1,$$
since $\ell_t-2q \le s-2q\le (\ell_j-1)-q$.
This means that among vertices of $C$ we can find a colour $i$ avoiding path of order $\ell_i-2q$ for some $i \neq j$ (note that technically we merge two colours in order to obtain a $t-1$ colouring).

Observe that there are at least $s-q-(\ell_i-2q)\ge q$ unused vertices left in $C$ so once again we can extend our path by alternating between $q$ vertices in $Q$ and these $q$ vertices of $C$. Since all the edges between $Q$ and $C$ are in colour $j\neq i$ the new path is still colour $i$ avoiding and has order $\ell_i$ as desired. 

\noindent Let us take a vertex $x$ and colour $i$ such that the degree of $x$ in the colour $i$ is maximum among all choices of $x$ and $i$. 

\begin{clm}\label{lemmaxdegree}
Suppose there is a path in $G-x$ of order $\ell_j-2$ avoiding some colour $j\neq i$. Then there is a colour $j'$ avoiding path of order $\ell_{j'}$ for some $j'$.
\end{clm}
\noindent{\it{Proof of Claim \ref{lemmaxdegree}.}}
Let $P$ be the colour $j$ avoiding path of order $\ell_j-2$ contained in $G-x$. Suppose $y,z$ are the endpoints of $P$. If both $xy,xz$ are not of colour $j$ then we can use $x$ to extend $P$ into a colour $j$ avoiding cycle of length $\ell_j-1$ and we are done by Claim \ref{lemcycle}. We are left with the following two cases.

\noindent{\bf  Case $1$.} Exactly one of the edges $xy,xz$ are of colour $j$, say w.l.o.g. $xz$.

If $x$ had any neighbours of colour $i$ outside $P$ then we can extend $zPyx$\footnote{Here by $aPb$ we denote the part of $P$ between two of its vertices $a$ and $b$.} to this vertex to find a colour $j$ avoiding path of order $\ell_j$ (since $i \neq j$). So we are done unless all colour $i$ neighbours of $x$ are on $P$. Let $w \neq y$ be such a neighbour. We claim that the neighbour $w^+$ of $w$ on $P$ closer to $y$ needs to be joined to $z$ by an edge of colour $j$. Indeed otherwise $xyPw^+zPwx$ is a colour $j$ avoiding cycle of length $\ell_j-1$ and we are done by Claim \ref{lemcycle}. This means that for any neighbour in colour $i$ of $x$ we know $z$ has a corresponding neighbour in colour $j$. In particular, within $P$ we know $z$ has at most one neighbour in colour $j$ less than $x$ has in colour $i$ (since $xy$ might be of colour $i$). In addition $x$ is a neighbour of $z$ in colour $j$ so within $P \cup \{x\}$ the colour $j$ degree of $z$ is at least as big as that of $x$ in colour $i$. Finally, take a vertex $a$ outside $P \cup \{x\}$ (note that it exists since $s \ge \ell_t >\ell_j-1=|P \cup \{x\}|$). If $za$ is not of colour $j$ then $xyPza$ is the desired path of order $\ell_j$ avoiding colour $j$ and we are done. Otherwise $z$ has more neighbours in colour $j$ than $x$ in colour $i$ (recall that $x$ has no colour $i$ neighbours outside $P \cup \{x\}$), which is a contradiction to the maximality of our choice of $x$ and $i$. 

\noindent{\bf  Case $2$.} Both edges $xy$ and $xz$ are of colour $j$.

Suppose $w$ is a neighbour of colour $i$ of $x$ on $P$. Note that $w \neq y$ since $xy$ is of colour $j\neq i$. As in the previous case consider the neighbour $w^+$ of $w$ on $P$ closer to $y$. If $zw^+$ is not of colour $j$ then we can rotate and find a path $yPw^+zPw$ of same length as $P$, which still avoids colour $j$ but of its endpoints $y,w$, exactly one is adjacent to $x$ in colour $j$ (since $xw$ is of colour $i \neq j$). So we can reduce to the first case and are done. This means that $zw^+$ is of colour $j$. This means that $z$ has at least as many neighbours of colour $j$ on $P$ as $x$ has neighbours of colour $i$.  

Note also that for any vertex $a$ outside of $P \cup \{x\}$ if $xa$ is of colour $i$ then $za$ must be of colour $j$, as otherwise $yPzax$ is the desired path avoiding colour $j$. So $z$ has at least as many neighbours of colour $j$ outside $P \cup \{x\}$ as $x$ has neighbours of colour $i$. Finally, since $x$ itself is a neighbour of $z$ in colour $j$ once again $z$ has more neighbours in colour $j$ than $x$ has in colour $i$, contradicting maximality of $x$ and $i$.  $\hfill \dashv$

\noindent We are now in a position to finish the proof. By Lemma \ref{lemi<3}, we may assume $\ell_i\ge 4$ for all $i$.

First by Claim \ref{removal} we have 
$$p({\ell_1-2}, \ldots, {\ell_{i-1}-2}, {\ell_i},  {\ell_{i+1}-2}, \ldots, {\ell_{t}-2}) \le \max\{\ell_{i},\ell_t-2,s-1\}$$

In case this is at most $s-1$ we can find either a colour $i$ avoiding path of order $\ell_i$ and are done or a colour $j$ avoiding path of order $\ell_j-2$ in $G-x$ and are done by Claim \ref{lemmaxdegree}. Otherwise, since $\ell_i\le \ell_t \le s,$ we must have $\ell_i=\ell_t=s$. This implies $s=(\ell_i-1)+1$ and we may apply Claim \ref{induction} with $j=i$ and $q=1$ to show  
$$p({\ell_1-2}, \ldots, {\ell_{i-1}-2},{\ell_{i+1}-2},\ldots  {\ell_{t}-2}) \le \max\{\ell_{t-1}-2,s-1\}\le s-1.$$
Once again this allows us to find a colour $j$ avoiding path of order $\ell_j-2$ for some $j \neq i$ inside $G-x$ and we are done by Claim \ref{lemmaxdegree}. This completes the proof.
$\hfill\Box$



\section{Concluding remarks and open problems}
In this paper, we determine the $(t-1)$-chromatic Ramsey number of a path of any length exactly, in other words we determine the minimum number of vertices $n$ needed for a $t$-edge coloured $K_n$ to contain a path of length $\ell$ using only $t-1$ of the colours. 

A natural further question could be to determine how the answer behaves if we want our path to use $s$ colours. 
\begin{qn}
Let $t>s \ge 1$ be integers. What is the smallest $n$ such that in any $t$-edge colouring of $K_n$ we can find a path of order $\ell$ using at most $s$ colours?
\end{qn}

When $s=1$ this recovers the classical Ramsey problem and the answer is known, up to lower order terms, to be between $(t-1)\ell$ and $(t-1/2)\ell$ and closing this gap is a major open problem in the area. We settle the opposite extreme, namely when $s=t-1$ and the answer is much smaller at $\ell+\floor{\frac{\ell-2}{2^t-2}}$. It seems like a nice question to determine how the answer transitions as $s$ varies from one extreme to another. 

A natural first step in attacking this question would be to try to resolve it if we seek a matching of size $\ell$ instead of $P_{2\ell}$. This question was raised by Gy\'arf\'as, S\'{a}rk\"{o}zy and Selkow \cite{gyarfas}.

A particularly nice instance of our result is that if $n \le 2^t$ then in any $t$ colouring of $K_n$ we can find a Hamilton path avoiding some colour. A nice question in this direction is whether the same holds for a Hamilton cycle?

In general it could be interesting to determine what happens in the $(t-1)$-chromatic case for cycles.
\begin{qn}
What is the smallest $n$ such that in any $t$-edge colouring of $K_n$ we can find a cycle of length $\ell$ using at most $t-1$ colours?
\end{qn}
Understanding the answer for other sparse structures, such as general trees might also be interesting.

\footnotesize

\bigskip

\vspace{0.5cm}

{\footnotesize
    
    \noindent Matija  Buci\'{c}
    
   \noindent Department of Mathematics, ETH, \"{Z}urich, Switzerland 
   
   \noindent \texttt{matija.bucic@math.ethz.ch}
    
   \medskip

   \noindent Amir Khamseh

    \noindent Department of Mathematics, Kharazmi University, 15719-14911 Tehran Iran
      
      
      \noindent \texttt{khamseh@khu.ac.ir}}

\end{document}